\documentclass{elsarticle}
\usepackage{amsfonts}
\usepackage[cp1250]{inputenc}
\usepackage{amsmath}
\usepackage{ifsym}
\usepackage{amsthm}

\newtheorem{thm}{Theorem}
\newtheorem{lemat}{Lemma}
\newtheorem{defi}{Definition}

\newtheorem{qu}{Question}

\author{Alan Czuro\'{n}
\\
Institute of Mathematics, Polish Academy of Sciences\\
\'{S}niadeckich 8\\
00-656 Warszawa, Poland\\
E-mail: alanczuron@gmail.com\\
+48 502-559-376}
\title{Property $F\ell_q$ implies property $F\ell_p$ for $1<p<q<\infty$}
\begin{document}

\begin{abstract}
It is known that for $\sigma$-compact groups Kazhdan's Property $(T)$ is equivalent to Serre's Property $FH$. Generalized versions of those properties, called properties $(TB)$ and $FB$, can be defined in terms of the isometric representations of a group on an arbitrary Banach space $B$. Property $FB$ implies $(TB)$.

It is known that a group with Property $(T\ell_p)$ shares some properties with Kazhdan's groups, for example compact generation and compact abelianization. Moreover in the case of discrete groups, Property $(T\ell_p)$ implies Lubotzky's Property $(\tau)$. 

 In this paper we prove that in the case of discrete groups and $\ell_p(\mathbb{N})$ spaces, for $1<p<q<\infty, p \not= 2$, Property $F\ell_q$ implies Property $F\ell_p$.
\end{abstract}
\begin{keyword} $F\ell_p$ property, Kazhdan's property $(T)$, Serre's fixed point property, affine isometric action, 1-cohomology.
\end{keyword}
\maketitle
\section{Introduction}

Property $(T)$ introduced by Kazhdan \cite{Ka} in terms of unitary representations, became a fundamental rigidity property of groups, with wide range of applications. It was proved by Delorme  \cite{De} and Guichardet \cite{Gu} that Kazhdan's Property $(T)$ is equivalent to Serre's Property $FH$ for $\sigma$-compact groups. Generalized versions of property $(T)$ and property $FH$, called properties $(TB)$ and $FB$, were introduced in \cite{BFGM} by Bader, Furman, Gelander and Monod, in terms of isometric representations of a group on an arbitrary Banach space $B$. Groups with those properties share some important properties with Kazhdan's groups, for example groups with property $(T{\ell_p})$ are compactly generated and have compact abelianization \cite[Theorem~6]{BO}. 

In this article, we study property $F\ell_p$, a fixed point property for affine actions on a real $\ell_p=\ell_p(\mathbb{N})$ spaces. 

Our main result is the following.

\begin{thm}
For every discrete countable group and $1<p<q<\infty, p\not=2$, property $F{\ell_q}$ implies $F{\ell_p}$.
\end{thm}

Fixed points properties for groups acting on general $L_p$ spaces for $p>2$ are poorly investigated. In this article we focus on the special case of a group actions on $\ell_p(\mathbb{N})$ spaces. At the begining we would like to present the current knowledge of the general property $FL_p$.

In \cite{BFGM} authors proved that higher rank algebraic groups and their lattices have fixed points for every affine action on $L_p$ spaces for $p>1$. Mimura \cite{Mi} showed that $SL_n(\mathbb{Z}[x_1,x_2,...,x_d])$ groups have fixed points for every affine isometric action on $L_p$ for $p>1$ and $d\ge 4$. It is known \cite{BFGM} that for every Kazhdan's group $G$ there exists a constant $\epsilon(G)>0$ such that every affine isometric action on $L_p$ spaces has a fixed point for $p\in [2,2+\epsilon(G))$. It is also known \cite{NS} that Gromov's random groups, containing $p$-expanders in their Cayley graphs have a fixed point property for affine actions on $L_p$ spaces for $p>1$. In \cite{N2} Nowak obtained sufficient conditions in terms of $p$-Poincar\'{e} constants implying that every affine isometric action of a given group on a reflexive Banach space has a fixed point. We refer to \cite{N1} for a recent survey.

It is known \cite{CDH} that the fixed point property for affine isometric group actions on $L_p(X,\mu)$ space, where $X$ is a space with measured walls implies property $(T)$. In the paper \cite{CDH} authors stated the following question.

\begin{qu}
Is the set of values of $p\in (1,\infty)$ for which a group has property $FL_p$ an interval?
\end{qu}

Theorem 1 partially answers this question in the case of discrete groups and $\ell_p$ spaces. 
\\
{\bf Acknowledgements.} I want to thank Piotr Nowak and Micha\l \hspace{1pt} Wojciechowski for valuable help during the work on this paper. This paper is supported by the WCNM.

\section{Preliminaries}

We recall the definition \cite{BFGM} of property $FB$ in the special case of discrete groups and $\ell_p$ spaces. By $\ell_p=\ell_p(\mathbb{N})$ we denote the usual Banach space of $p$-summable real sequences. For $v\in \ell_p(\mathbb{N})$ we use the following notation $v=(v[i])_{i=1}^{\infty}$.

\subsection{Representations and cocycles}
Denote by $O(\ell_p)$ the group of linear bijective isometries of the Banach space $\ell_p$. Let $G$ be a discrete group.
\begin{defi}
(Isometric representation) An isometric representation of $G$ on Banach space $\ell_p$ is homomorphism
$$\pi:G \rightarrow O(\ell_p).$$
\end{defi}
We denote the image of $g\in G$ under the homomorphism $\pi$ by $\pi (g)$.
The affine group $\mathrm{Aff}(\ell_p)$ of a real affine space (a vector space who forgot its origin) consists of invertible maps satisfying

$$T(tv+ (1-t)w)=tT(v)+ (1-t)T(w),$$
for $t\in \mathbb{R}$ and $v,w\in \ell_p$.

An isometric affine action of a group $G$ on $\ell_p$ is a homomorphism $G\rightarrow \mathrm{Aff}(\ell_p)$ of the form:
$$g(v)=\pi(g)(v) + b(g),$$
where $\pi$ is an isometric representation and $b:G\rightarrow \ell_p$ is a $\pi$-cocycle, that is, an element of the Abelian group
$$Z^{1}(\pi)=\{ c:G\rightarrow \ell_p: c(gh)=\pi (g) (c(h))+c(g); g,h \in G \}.$$

\begin{defi}
We say that a group $G$ has the property $F\ell_p$ if any continuous action of a group $G$ on $\ell_p$ space by affine isometries has a $G$-fixed point.
\end{defi}

Group $Z^1(\pi)$ contains the subgroup of $\pi$-coboundaries:

$$B^1(\pi)= \{ b(g)=\pi (g)(v)-v: v\in \ell_p \}.$$

Observe that $Z^1(\pi)$ describes all affine actions with linear part $\pi$, and $B^1(\pi)$ corresponds to those actions which have a $G$-fixed point. This interpretation involves the choice of origin point in the space. Two cocycles differing by a coboundary can be thought of defining the same affine action viewed from different reference points. 

To prove that a group $G$ has the property $F\ell_p$ it is sufficient to show that for every isometric representation $\pi:G\rightarrow O(\ell_p)$, we have
 $$H^1(G,\pi)=Z^1(\pi)/B^1({\pi})=\{ 0\}.$$
The group $H^1(G,\pi)$ is the first cohomology group of $G$ with $\pi$-coefficients.
In this spirit we can reformulate our result as follows.

\begin{thm}
The set of parameters $p\in (1,\infty)$ for which $H^1(G, \pi)=\{0 \}$ for every representation $\pi: G \rightarrow O(\ell_p)$ is an interval or an interval without the point $\{ 2\}$. 
\end{thm}

In our proof we will use the following results.

\begin{thm}
(Banach-Lamperti) Assume that $p \neq 2$. Every linear isometry $\pi:\ell_p \rightarrow \ell_p$ has the following form
$$\pi(v[i])=\epsilon(i)v[\sigma(i)],$$
where $\sigma:\mathbb{N}\rightarrow \mathbb{N}$ is bijection of natural numbers and $\epsilon:\mathbb{N}\rightarrow \{ -1, 1  \}$.
\end{thm}

Let $p\in (1,\infty)\backslash \{2 \}$ and $\pi :G\rightarrow O(\ell_p)$ be an isometric representation. By the Banach-Lamperti theorem for every $g\in G$ the operator $\pi (g)$ is given by the formula

\begin{equation}\label{mapa}
\pi(g)(v[i])=\epsilon_g(i) v(\sigma_g(i)).
\end{equation} 
This formula defines a representation on an arbitrary $\ell_p$ space, for $1\le p \le \infty$. We denote this representation by $\pi$, no matter on which $\ell_p$ the representation is considered.
 
\subsection{Graphs and Inequalities}
In this section we consider finite connected graphs. Denote by $V$ the set of vertices and by $E$ the set of edges.
\begin{defi}
Let $X=(V,E)$ be a graph. Given  $A\subset V$ define its edge boundary $\partial^e (A)$ setting $\partial^e (A) =\{ e \in E:\textit{e has exactly one vertex in A}   \}$
\end{defi}

Now we define the so-called Cheeger constant.

\begin{defi}
The Cheeger constant $h(X)$ of a finite graph $X$ is given by the formula:
$$h(X)= \min \{  \frac{\# \partial^e A}{\#A} : 0 < \# A\le \frac{\#X}{2}\}.$$ 
\end{defi}

\begin{defi}
Let $\{ X_n \}_{n=1}^{\infty}$, $X_n=(V_n,E_n)$ be a countable collection of finite graphs of uniformly bounded degree. Assume that there exists a constant $c>0$ such that
$$h(X_n)>c,$$
for every $n\in \mathbb{N}$. Then the collection $\{ X_n \}$ is said to be a uniformly expanding family of graphs.
\end{defi}

In this paper we consider graphs as metric spaces. The distance between two vertices is the number of edges in a shortest path connecting them. Let $B(x,r)$ denote the ball of radius $r$ centered at a point $x$.
For a uniformly expanding family of graphs the following Poincar\'{e} type inequality holds.

\begin{thm}
(J. Matou\v{s}ek) Let $X_n=(V_n,E_n)$ be a uniformly expanding family of graphs of bounded degree and $p\in (1, \infty)$. Then there exists $c(p)>0$ such that for every function $f:V_n\rightarrow \mathbb{R}$ with average 0,
$$ \frac{  \sum_{x\in V_n} \sum_{y\in B(x,1)}|f(x)-f(y)|^p }{  \sum_{x\in V} |f(x)|^p  }  > c(p). $$
\end{thm}

This inequality is a generalisation of Alon-Milman inequality, see  \cite{Do}, \cite{Al}, \cite{Tan} and \cite{AM} for more details. Theorem 4 was proved by Matou\v{s}ek in \cite[Proposition~3]{Mat}. In the proof author assumes that the collection of graphs is a d-regular expander. In fact it sufficies to assume that considered collection of graphs a uniformly expanding family of graphs. The proof of Matou\v{s}ek theorem in the case of uniformly expanding family of graphs is almost the same as in the case of expander graphs.

\section{Proof of the Theorem 1}
The rest of the paper is concerned with proving Theorem 1. 
\subsection{Fixed point of representation}
In this section we prove the following theorem.
\begin{thm}
Assume that the group $G$ has property $F\ell_q$, but fails to have property $F\ell_p$ for $1<p<q<\infty$, $p\not=2$. Let $\pi$ be a representation on $\ell_p$ space which admits a non-trivial cocycle. Let $r=\frac{q^2}{p}$. Then there exists a vector $z\in \ell_r$ which is a fixed point of representation $\pi$ on $\ell_r$ space and $z\not\in \ell_q$.
\end{thm}
We start by some simple observations. Let $p\in (1,\infty)\backslash \{2\}$ and $p<q$. Let $b^p:G\rightarrow \ell_p$ be a cocycle associated to the representation $\pi$ on $\ell_p$ space. Let $i:\ell_p \hookrightarrow \ell_q$ be the cannonical inclusion. It is easy to check that $b^q=i(b^p)$ is a cocycle associated to the representation $\pi$ on $\ell_q$.

On the other hand let $b^q$ be a cocycle associated to the representation $\pi$ on $\ell_q$ space. Assume that for every $g\in G$ we have
$$||b^q (g)||_p < \infty.$$
Then $b^q:G\rightarrow \ell_p$ is a cocycle associated to the representation $\pi$ on $\ell_p$ space. Now we pass to the proof of theorem 5.

\begin{proof}

Assume that there exist $p,q \in \mathbb{R}$ satisfying $1<p<q<\infty$, $p\not=2$ such that $G$ has property $F\ell_q$ but fails to have property $F\ell_p$.  Since $G$ does not have Property $F\ell_p$ we know that there exists a representation $\pi$ on $\ell_p$  which admits a non-trivial cocycle
$$b^p:G\rightarrow l_p,$$
which by definition implies that:
\begin{lemat}
For every $v\in \ell_q$ satisfying $b^p(g) =\pi(g)(v)-v$ where $g\in G$ we have
 $$v \not\in \ell_p.$$
\end{lemat}

Since $G$ has property $F\ell_q$ and $b^p$ is a cocycle associated to the representation $\pi$ on $\ell_q$, there exists a vector $v\in \ell_q$ such that
\begin{equation}\label{2}
b^p(g)=\pi(g)(v)-v.
\end{equation}

This implies that for every $g\in G$, we have $\pi(g)(v)-v=b^p(g) \in \ell_p$, but by Lemma 1 we have $v\not\in \ell_p$. Hence denoting $v=(v[i])_{i=1}^{\infty}$, and $\pi(g)(v)=( \epsilon_g(1) v[\sigma_g(1)],\epsilon_g(2) v[\sigma_g(2)],...)$ we have:
\begin{equation}\label{3}
\sum_{i=1}^{\infty} |v[i]|^p = \infty,
\end{equation}
and 
\begin{equation}\label{4}
\sum_{i=1}^{\infty} |\epsilon_g(i) v[\sigma_g(i)] -v[i]|^p < \infty.
\end{equation}

Now we show that \eqref{3} and \eqref{4} lead to existence of "non-trivial" fixed point of the representation $\pi$ on $\ell_r$ space. In the proof of this fact we can assume that for every  $i\in \mathbb{N}$ and $g\in G$ we have $v[i] \ge 0$ and $\epsilon_g (i)=1$. Indeed, since:
$$
\sum_{i=1}^{\infty} |\epsilon_g(i) v[\sigma_g (i)]-v_i|^p \ge \sum_{i=1}^{\infty} \Bigg| |\epsilon_g(i) v[\sigma_g (i)]|-|v[i]| \Bigg|^p,
$$
our assumption does not change $\eqref{3}$ and $\eqref{4}$. Put $w=(w[i])_{i=1}^{\infty}$ where $w[i]=(v[i])^{\frac{p}{q}}$.
Obviously $w\not\in \ell_q$, but $w \in \ell_{r}$ where $r=\frac{q^2}{p}$ and $\pi(g)(w)- w\in \ell_q$, for every $g\in G$. Indeed, since for $\alpha<1$
$$||x|^{\alpha}-|y|^{\alpha}|\le ||x|-|y||^{\alpha},$$ 
we get:

 $$||\pi(g)(w)- w||_{q}^{q}=\sum_{i=1}^{\infty} |w[\sigma_g(i)] - w[i]|^q =\sum_{i=1}^{\infty} |v[\sigma_g(i)]^{\frac{p}{q}} - v[i]^{\frac{p}{q}}|^q$$
 $$\sum_{i=1}^{\infty} |v[\sigma_g(i)]^{\frac{p}{q}} - v[i]^{\frac{p}{q}}|^q 
  \le \sum_{i=1}^{\infty} (|v[\sigma_g(i)]-v[i]|^{\frac{p}{q}})^q =\sum_{i=1}^{\infty} |v[\sigma_g(i)]-v[i]|^p < \infty .$$
 
Now we use the assumption that $G$ has property $F\ell_q$. Since $\pi(g)(w)-w\in \ell_q$, it defines a cocycle $c$ in $\ell_q$. Since $G$ has Property $F\ell_q$ it follows that there exists $u \in \ell_q$ such that 
\begin{equation}
\pi(g)(w)-w=c=\pi(g)(u)-u,
\end{equation}
or equivalently
\begin{equation}
\pi(g)(w-u)=w-u.
\end{equation}
i.e. $z=w-u$ is a fixed point of a representation $\pi$ on $\ell_r$ space. Observe that $z$ belongs to $\ell_r$, but $z\not\in \ell_q$, because $w\not\in \ell_q$. This ends the proof of the Theorem 5.
\end{proof}

The assumption that $v[i]\ge 0$ and $\epsilon_g(i)=1$ is purely technical and not necessary. If we omit this assumption than vector $w$ should have the following form $w[i]=sgn(v[i])|v[i]|^{\frac{p}{q}}$.

\subsection{Construction of Components}

In this section we use Theorem 5 to decompose the natural numbers $\mathbb{N}$ onto disjoint \emph{components} closed under permutations $\sigma_g$ associated to the representation $\pi$.
Assume that the representation $\pi$ and the vectors $v$, $w$ and $z$ are as in the proof of Theorem 5. 

Let $z=(z[i])_{i=1}^{\infty}$. Put $J=\{ i: z[i]=0 \}$ and let $\mathbb{N}-J=\cup_{I \in P} I $  be the decomposition into disjoint orbits of the representation $\pi$, where $P$ is the (countable) family of \emph{components}. Here by \emph{component} we mean any smallest by inclusion set closed under permutations $\sigma$ induced by the representation. Note that all components $I\in P$ are finite. Indeed, since $z$ is a fixed point of the representation $\pi$ on $\ell_r$ space, it follows that $z[i]\not=0$ is constant on every $I \in P$. But $z\in \ell_r$, every component has to be finite. On the other hand, $P$ is infinite because otherwise $z$ would have finite support, which contradicts the fact $z\not\in \ell_q$.

From now on by \emph{components} we mean disjoint orbits of permutations $\sigma$ associated to the representation $\pi$ constructed as above. 

\subsection{Construction of graphs and proof of Theorem 1}
In this section we construct a family of graphs associated to the representation $\pi$. The assumptions on $\pi$ are as in Section $3.1$. Let $P$ be the family of components constructed in Section $3.2$.

We have assumed that the group $G$ has Property $F\ell_q$. This implies \cite[Theorem~6]{BO} that $G$ is finitely generated. Let $S$ be a symmetric set of generators of the group $G$. Consider the collection of finite connected graphs $X_I = (V_I, E_I)$ for $I\in P$, where $V_I=I$ and vertices $a,b\in V_I$ are connected by an edge if there exists a generator $g\in S$  such that $\pi(g)({\bf 1}_a)={\bf 1}_b$    i.e., for some $g\in S$ the permutation $\sigma_g$, associated to the representation $\pi$ satisfies $\sigma_g(a)= b$. Note that two vertices may be connected by more than one edge. Observe that the collection $\{X_I\}_{I\in P}$ of graphs has uniformly bounded degree.

Observe that to prove Theorem 1 it sufficies to prove two following theorems.

\begin{thm}
Assume that the collection $\{X_I\}_{I\in P}$ is not a uniformly expanding family of graphs. Then the representation $\pi$ on $\ell_q$ space admits a non-trivial cocycle, for every $1\le q< \infty$.
\end{thm}

\begin{thm}
Let $G$ have property $F\ell_q$. Assume that the collection $\{X_I\}_{I\in P}$ associated to the representation $\pi$ on $\ell_q$ space is a uniformly expanding family of graphs. Then every cocycle associated to the representation $\pi$ on $\ell_p$ space for $1\le p \le q$ is a coboundary.
\end{thm}

The proof of the Thorem 1, assuming Theorems 6 and 7 is given in the next section.

\subsection{Proof of the Theorem 1}

Let $1\le p<q <\infty$ and $(p\not=2)$. Assume by contrary that a group $G$ has Property $F\ell_q$ but fails to have Property $F\ell_p$. Let $\pi$ be the representation on $\ell_p$ that admits a non-trivial cocycle $b^p$. Let $\{ X_I\}_{I\in P}$ be collection of graphs constructed as in Section $3.3$. Assume that this collection is not a uniformly expanding family of graphs. The by Theorem 6 the representation $\pi$ on $\ell_q$ admits a non-trivial cocycle, thus $G$ does not have property $F\ell_q$. This implies that under our assumptions constructed collection of graphs is a uniformly expanding family of graphs. In this case Theorem 7 implies that every cocycle associated to the representation $\pi$ is a coboundary. This contradicts our assumption that $\pi$ admits a non-trivial cocycle on $\ell_p$. This ends the proof of the Theorem 1.
\subsection{Proof of the Theorem 6}

Assume that $\pi$ is the representation of the group $G$ as in Section $3.1$. Let $\{X_I \}_{I\in P}$ be a collection of graphs constructed as in Section $3.3$. Assume that the collection $\{X_I \}_{I\in P}$ is not a uniformly expanding family of graphs. We construct a non-trivial cocycle for the representation $\pi$ on every $\ell_q$ space, for $1\le q < \infty$.

The collection $\{X_I \}_{I\in P}$ of graphs of uniformly bounded degree is not a uniformly expanding family of graphs, thus for every $\delta>0$ there exists $I\in P$ such that
$$h(X_I)<\delta.$$
On the other hand every graph $X_I$ is connected, thus 
$$h(X_I)>0,$$
for every $I\in P$.
This implies that there exists an infinite family $P' \subset P$ and subsets $A_I\subset I$ for every $I\in P'$, such that $\#A_I \le \frac{\#I}{2}$ and  $\sum_{I\in P'} \frac{\# \partial^e A_I}{\# 
A_I}<\infty$.

Consider $v=(v[i])_{i=1}^{\infty} \in \ell_{\infty}$ where $v[i] =(\frac{1}{\#A_{I}})^{\frac{1}{q}}$ for $i\in A_I$ and $v[i]=0$ otherwise. Observe that $v\not\in \ell_q$. 
 For every generator $g\in S$ define the sets $\Omega[g], \Omega[g]'\subset \mathbb{N}$ as follows. 
 
 Index $i\in \Omega[g]$ if there exists $I\in P'$ such that $i \in A_I$ and $\sigma_g(i)\not \in A_I$.
 
 Index $i\in \Omega'[g]$ if there exists $I\in P'$ such that $i \not\in A_I$ and $\sigma_g(i) \in A_I$.

Recall that $S$ is a symmetric set of generators of the group $G$. For every $g\in S$ we have
 
$$
\aligned
||\pi(g)(v)-v||_q^q =& \sum_{i=1}^{\infty} |v[\sigma_g(i)]-v[i]|^q\\=&
 \sum_{I\in P'} \sum_{i\in I} \big( \sum_{\Omega[g]} (\frac{1}{\# A_{I}})+ \sum_{\Omega'[g]} (\frac{1}{\# A_{I}}) \big)\\ \le&
2\sum_{I\in P'} \frac{\partial^{e} A_{I}}{\# A_{I}} < \infty
\endaligned
$$
 
Thus  $\pi(g)(v)-v\in \ell_q$ for $g\in S$. Now we show that $\pi(g)(v)-v\in \ell_q$ for every $g\in G$. For every $g\in G$ we have $g=s_1s_2...s_l$ for $s_i \in S$. We can consider $s_1s_2...s_l$ as a finite path connecting neutral element $e\in G$ with $g$ in the Cayley graph of $G$ with respect to the set of generators $S$. Put $s_0 = e$. We have
$$
\pi(g)(v) - v = \sum_{k=0}^{l-1} \pi (s_0s_1...s_{k+1})(v) - \pi(s_0s_1...s_{k})(v) = 
\sum_{k=0}^{l-1} \pi(s_0s_1...s_k)(\pi(s_{k+1})(v) - v)
$$

Since $\pi$ is an isometric representation we get $\pi(g)(v)-v\in \ell_q$ for $g\in G$. It follows that $b(g)^q=\pi(g)(v)-v$ is a cocycle with respect to the representation $\pi$ on $\ell_q$ space.
We have to show that the cocycle $b^q(g)=\pi(g)(v) - v$ is not a coboundary on $\ell_q$.

 Assume, by contrary, that $\pi(g)(v)-v$ is a coboundary. Then there exists a vector $w\in \ell_q$ such that 
 $$\pi(g)(w)-w=\pi(g)(v)-v,$$
 for every $g\in G$ i.e. $w-v$ is fixed point of representation $\pi$. Thus there exists $c_I, d_I \in \mathbb{R}$, such that the vector $w$ has the following form

\begin{equation}
w_i=\begin{cases} (\frac{1}{A_{I}})^{\frac{1}{q}} + c_{I} &\text{for } i\in A_{I} \\
c_{I} &\text{for }  i\in V_{I}- A_{I}, I\in P'\\
d_I &\text{for } i \in V_{I}, I\in P\backslash P'
\end{cases}
\end{equation}

It is clear that
$$||w||_{q}^{q}=\sum_{I\in P}||w_{|_{I}}||_{q}^{q}.$$

For every $I\in P'$ either 
 $c_I\ge -(\frac{1}{\# A_I})^{\frac{1}{q}}/2$ 
 or
 $c_I\le -(\frac{1}{\# A_I})^{\frac{1}{q}}/2$.
 In the first case we have for every $I\in P'$ we have
\begin{equation}
\aligned
||w_{|_{I}}||_q^q = \sum_{i\in I}|w[i]|^q &= \sum_{i\in A_I}|w[i]|^q +\sum_{i\in \{V_I-A_I\}}|w[i]|^q\\&\ge
  \sum_{i\in A_I}|\frac{1}{2}(\frac{1}{\#A_I})^{\frac{1}{q}}|^q + \sum_{i\in \{ V_I-A_I\}}|w[i]|^q \\&\ge
   \frac{1}{2^q}\sum_{i\in A_I}\frac{1}{\#A_I}+\sum_{i\in \{ V_I-A_I\}}|w[i]|^q\\&\ge 
 \frac{1}{2^q}||v_{|_I}||_q^q.
\endaligned
\end{equation}
 Similarly in the second case using the fact that $\#\{V_I-A_I \}\ge \# A_I$ we have

\begin{equation}
\aligned
||w_{|_{I}}||_q^q = \sum_{i\in I}|w[i]|^q &= \sum_{i\in A_I}|w[i]|^q +\sum_{i\in \{V_I-A_I\}}|w[i]|^q\\&\ge
   \sum_{i\in A_I}|w[i]|^q+\sum_{i\in \{V_I-A_I\}}|\frac{1}{2}(\frac{1}{\#A_I})^{\frac{1}{q}}|^q  \\&\ge 
   \sum_{i\in A_I}|w[i]|^q+\frac{1}{2^q}\sum_{i\in A_I}\frac{1}{\#A_I}\\&\ge 
 \frac{1}{2^q}||v_{|_I}||_q^q.
\endaligned
\end{equation}

Thus there does not exists vector $w\in \ell_q$ such that $b^q(g)=\pi(g)(w)-w$. This prove that $w\not\in \ell_q$. Thus cocycle $\pi(g)(v)-v$ is not a coboundary, which ends the proof.

\subsection{Proof of the Theorem 7}
Assume that $\pi$ is the representation of a group $G$ as in Section $3.1$. Let $P$ be the decoposition of natural numbers onto disjoint components as in Section $3.2$. Let $\{X_I \}_{I\in P}$ be collection of graphs constructed as in Section $3.3$. Assume that the collection $\{X_I \}_{I\in P}$ is a uniformly expanding family of graphs. We prove that every cocycle associated to the representation $\pi$ on $\ell_p$, for $1\le p \le q$ is a coboundary.

Let $1\le p <q$. Assume that $b^p$ is a non-trivial cocycle associated to the representation $\pi$ on $\ell_p$ space. Let us recall that Property $F\ell_q$ implies that there exists $v\in \ell_q$ such that:
$$b^p(g)=\pi(g)(v)-v.$$ 
 
Consider the vector $\widehat{v}=(\widehat{v}[i])_{i=1}^{\infty}$, defined by $\widehat{v}[i] = v[i] - \frac{1}{\# I}\sum_{i\in I} v[i]$, for $i\in I$, where $I\in P$. It follows from the construction of $P$ that 
  $$\pi(g)(v)-v=\pi(g) (\widehat{v})-\widehat{v}.$$
Observe that we can treat $\widehat{v}_{|_I}$ as a function on the graph $X_I$. Moreover, this function has average zero. Thus by Theorem 4 we have:
$$
\frac{  \sum_{x\in V_I} \sum_{y\in B(x,1)}|\widehat{v}[x]-\widehat{v}[y]|^p }{  \sum_{x\in V_I} |\widehat{v}[x]|^p  }  > c(p). 
$$  
Observe that 
$$
\sum_{g\in S} ||\pi(g) (\widehat{v})-\widehat{v}||_p^p=\sum_{x\in V_I} \sum_{y\in B(x,1)}|\widehat{v}[x]-\widehat{v}[y]|^p ,
$$
  thus, using the fact that $\pi(g) (\widehat{v})-\widehat{v}=b^p(g)\in \ell_p$, we have
$$
 \infty >\sum_{g\in S} ||\pi(g)(v)-v||_p^{p}= \sum_{g\in S} ||\pi(g) (\widehat{v})-\widehat{v}||_p^p \ge c||\widehat{v}||_p^p.$$
 Thus  $\widehat{v}\in \ell_p$ and $b^p=\pi(\widehat{v})-\widehat{v}$ - which contradicts the  assumption that $b^p$ is a non-trivial cocycle.

\end{document}